\documentclass{article}

\usepackage{arxiv}

\usepackage[utf8]{inputenc} 
\usepackage[T1]{fontenc}    
\usepackage{hyperref}       
\usepackage{url}            
\usepackage{booktabs}       
\usepackage{amsfonts}       
\usepackage{nicefrac}       
\usepackage{microtype}      
\usepackage{lipsum}
\usepackage{graphicx}
\graphicspath{ {./images/} }

\title{Sufficient  Loading Conditions for Self-organization in the $n$-dimensional Version of BML Model with Stochastic Direction Choice
}

\author{
 Valery V. Kozlov \\
Steklov 
 Mathematical Institute\\
of Russian Academy of Sciences\\
Moscow, Russial\\ 
  \texttt{kozlov@pran.ru } \\
 \And
 Alexander G. Tatashev  \\
  Department of Higher Mathematics\\
  Moscow Automobile and Road Construction\\
  State Technical University (MADI) \\
  Moscow, Leningradsky avenue, 64, Russia  \\
  \texttt{a-tatashev@yandex.ru} \\
   \And
 Marina V. Yashina  \\
  Department of Higher Mathematics\\
  Moscow Automobile and Road Construction\\
  State Technical University (MADI) \\
  Moscow, Leningradsky avenue, 64, Russia  \\
  \texttt{mv.yashina@madi.ru} \\
}

\begin{document}
\maketitle
\begin{abstract}
 
A dynamical system is considered, which comprises an $n$-dimensional lattice $N_1 \times N_2 \times \dots \times N_n$ with periodic boundary conditions. Particles traverse this lattice following a variant of the Biham--Middleton--Levine (BML) traffic model's particle movement rules. We have proved that the BML model, when treated as a dynamical system, constitutes a specialized class of a Buslaev net. 
This equivalence allows us to employ established Buslaev net analysis techniques to investigate the BML model.
In Buslaev nets conception the self-organization property of the system corresponds to the existance of velocity single point spectrum equal to 1.
One notable aspect of the model under consideration is that particles can change their type with a certain probability. For simplicity, we assume a constant probability $q$ that a particle changes type at each step.
In the l case where $q=0$, the system corresponds to the classical version of the BML model. We define a state of the system where all particles continue to move indefinitely, both in the present and the future, as a state of free movement. 
A sufficient condition for the system to result in a state of free movement from any initial state (condition for self-organization)  has been found.
This condition is that the number of particles be not greater than half the greatest common divisor of the numbers $N_1$ $N_2,\dots,N_n.$

\end{abstract}


\section{Introduction}
Urban traffic congestion is a critical issue affecting cities worldwide. It results in significant economic losses, environmental damage, and substantial disruptions in individuals' daily lives. This problem is particularly acute in large megacities, prompting extensive research efforts aimed at improving transportation infrastructure and minimizing public expenditures. The adoption of new automotive technologies plays a pivotal role in enhancing traffic control, while the development of communication networks offers drivers opportunities to select optimal routes, thereby alleviating congestion in major cities. Consequently, the development of algorithms for effective traffic control is of utmost importance. Mathematical models serve as invaluable tools for deriving optimal solutions to these challenges. Therefore, predicting and mitigating traffic jams are crucial focal points.
One widely utilized approach for modeling traffic flows in urban networks and predicting traffic congestion is based on various versions of the renowned Biham--Middleton-Levine (BML) model, ~\cite{Biham}. The BML model, renowned for its simplicity and efficiency, was introduced to address traffic-related issues.

 In O. Biham, A.A. Middleton and D. Levine study~\cite{Biham}, a mathematical traffic model, the BML model, was introduced. In this model, particles navigate a two-dimensional toroidal lattice comprising $N_1$ rows and $N_2$ columns, moving at discrete moments. 
These particles come in two types: particles of the first type (referred to as "red particles") move horizontally from left to right, while particles of the second type (referred to as "blue particles") move vertically along columns. At any given moment, a cell can either be vacant or contain a particle of one of these two types. Particle movement along rows or columns adheres to a rule derived from the elementary cellular automaton 184 (ECA 184)~\cite{Wolfram}. According to this rule, a particle advances to an adjacent vacant cell in the same direction, or it remains stationary if the cell ahead is occupied. Additionally, a rule is in place to prevent two particles from entering the same cell simultaneously, one moving along a row and the other along a column. Specifically, the rule dictates that particles of the first type move first, followed by particles of the second type. Alternatively, one can envision a scenario where particles of the first type move during odd time steps, and those of the second type move during even time steps.

In a study by Linesh et al~\cite{Linesh}, results are presented that investigate the behavior of the BML model across various lattice dimensions and particle densities. This investigation relies on simulation modeling techniques.
One of the primary challenges in analyzing traffic model systems is determining the average particle velocity, which quantifies the ratio of particles moving during a step to the total number of particles. If, starting from a certain point in time, all particles move without delay, the particle velocity reaches a value of 1, indicating a state of free movement within the system. This inherent property of the system to transition into a state of free movement is called the as self-organization.

In the work by Austin et al~\cite{Austin}, it has been provted that if $N_1=N_2=N$ and the number of particles $m$ is less than $N/2$, the system consistently evolves into a state of free movement from any initial state. Similarly, Moradi~\cite{Moradi} established that when $m= N/2$, the system also attains a state of free movement from any initial state.

For a related one-dimensional system, where particles adhere to the rule ECA184 while traversing a circular lattice with $N$ cells, Blank~\cite{Blank} and Belitsky and Ferrary~\cite{Belitsky} provided evidence that if the number of particles does not exceed $N/2$, the system reliably transitions into a state of free movement from any starting point. Conversely, if there are more than $N/2$ particles, the average velocity of particles drops below 1. Additionally, they derived a formula for calculating the average velocity.

In ~\cite{Tatashev}, a system is considered such that this system contains two closed one-dimensional lattices called contours (two-contour system). The lengths of contours are different. The system belongs to a class of dynamical systems called Buslaev nets
~\cite {Kozlov_01}, ~\cite{Kozlov_02}. 
 In ~\cite{Tatashev}, it is supposed that there is a common cell of the contours, and particles move along each contour according to the rule of ECA~184. After passing a node, a particle passes to another contour with some probability. It has been proved that, if the number of particles is not greater than a half of the greatest common divisor of the numbers of cells in the contours, then the system results in a state of free movement from any initial state. Under the additional condition that the probability of transition from the first contour circuit to the second contour and the probability of transition from the second circuit to the first are not equal to~0, this sufficient condition for self-organization is also necessary. In proving these statements, known facts of the number theory relating to the theory of equations of two variables in integers are used.

Tatashev and Yashina~\cite{Tatashev} explored a system involving two closed one-dimensional lattices, referred to as contours, within a two-contour system. These contours had varying lengths and belonged to the class of dynamical systems known as Buslaev nets~\cite{Kozlov_01, Kozlov_02}. In this study, it was assumed that there was a common cell shared by the contours, and particles traversed each contour following the ECA~184 rule. The research established that when the number of particles did not exceed half of the greatest common divisor of the cell numbers in the contours, the system reliably reached a state of free movement from any initial state. Moreover, when an additional condition stipulating nonzero transition probabilities from the first contour circuit to the second, and vice versa, was imposed, this sufficient condition for self-organization was also proven to be necessary. The proofs of these assertions relied on well-established principles from number theory, specifically relating to equations involving two integer variables.

It is worth mentioning that in Wenbin et al's investigation~\cite{Wenbin}, a version of the BML model was considered where particles could pass through common cells of both rows and columns simultaneously. In this scenario, the system was also represented as a contour network, and certain cells were identified as nodes while others were not.
 
Additionally, Ding et al~\cite{Ding} and Yashina et al~\cite{Yashina} explored stochastic versions of the Biham-Middleton-Levine model.

In another study by Chau et al~\cite{Chau}, an $n$-dimensional version of the BML model was introduced within an $N_1 \times N_2\dots \times N_n$ lattice, incorporating periodic boundary conditions. In reference ~\cite{Yashina}, the model was further examined through simulations under the assumptions that $n=3$ and $N_1=N_2=N_3$.

In the works of Criado and Wan~\cite{Criado} and Gillman and Martinsson~\cite{Gillman}, various mathematical models with a two-dimensional structure are explored, along with their applications.

In this paper, we extend upon the system discussed in Biham-Middleton-Levine's work~\cite{Biham}. This system encompasses an $n$-dimensional lattice with dimensions $N_1 \times N_2 \times \dots \times N_n$. At discrete moments in time, each particle undergoes a type change with a certain probability. For simplicity, we assume a constant probability denoted as $0 \leq q < 1$, and the new direction is chosen with equal probability among available options.

When $n=2$, $N_1=N_2$, and $q=0$, the system aligns with the one analyzed in Austin et al's study~\cite{Austin}. It has been established that a sufficient condition for the system to evolve into a state of free movement from any initial state is that the number of particles must not exceed half of the greatest common divisor of the dimensions $N_1$, $N_2$, and so on, up to $N_n$. The proof of this condition draws upon methodologies employed in Austin et al~\cite{Austin}, Moradi et al~\cite{Moradi}, Tatashev and Yashina~\cite{Tatashev}, and incorporates principles from number theory~\cite{Andrews}.

Furthermore, we present a description of the system in terms of Buslaev nets, building upon the established embedding of the BML model into the Buslaev nets class, as demonstrated in Yashina et al's work~\cite{Yashina}. In this framework, the rows of lattice cells are referred to as contours, with each cell corresponding to a node shared by multiple contours. This allows us to apply methodologies developed for the analysis of contour networks to the study of the BML model. A crucial concept introduced here is the spectrum of the Buslaev net, encompassing sets of limit cycles (closed trajectories within the net's state space) and the corresponding velocities of particles moving within the contours.

\section{
Description of model 
}
\label{section:Desc}

The system contains an $n$-dimensional 
lattice $N_1\times N_2\times \dots \times N_n$ with periodic boundary conditions.

\begin{figure}[ht!]
\centerline{\includegraphics[width=200pt]{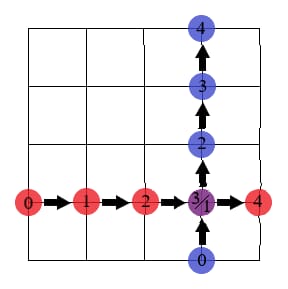}}
\caption{ Classical  BML model trajectory with $n=2$}
\end{figure}

 There are $m<N_1N_2\dots N_n$ particles. At each moment, any particle belongs to one of $n$ types. In general case, the type of a particle may be changed.  Particles move at discrete moments, $t=0,1,2,\dots$ At each moment of time,  first particles of the type~1 move, then particles of the type~2 move, then particles of the type~3 move, etc, Fig. 1.

After that, type of each particle changes with probability $q$, $0\le q< 1.$ Suppose, at time $t,$ a particle of the $i$th type is in the cell  $(x_1,\dots,x_n)$ and, after that the particles of  
the types $1,\dots,i-1$ could move, the cell $(x_1,\dots,x_{i-1},x_i+1,x_{i+1},\dots,  x_n)$ is vacant, then, at time $t+1,$ the particle will be in the cell  
$(x_1,\dots,x_{i-1},x_i+1,x_{i+1},\dots,  x_n),$ $i=1,\dots,n,$ $0\le x_j\le N_j,$
$j=1,\dots,N,$  , Fig. 2 with $n=2.$ 

Note that, if $n=2,$ then we have a BML model in two-dimensional torroidal lattice.

\begin{figure}[ht!]
\centerline{\includegraphics[width=300pt]{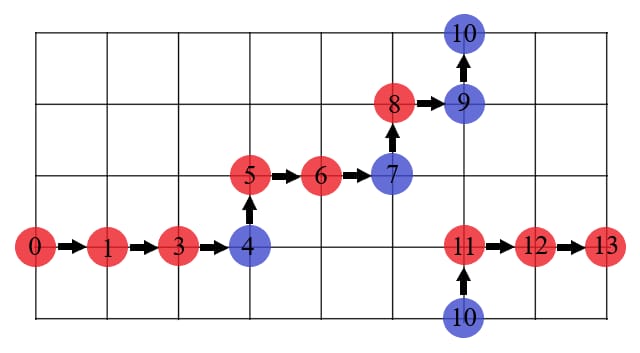}}
\caption{ Trajectories in BML model with $q$ probability of direction change}
\end{figure}

\section{
BML model in terms of Buslaev nets
}
\label{section:BuslN}

Suppose $q=0.$ Let us describe the model in terms of Buslaev nets.

There are $N_1\cdot N_2\cdot N_{i-1}\cdot \dots \cdot N_{i+1}\cdot N_n$ contours
called the contours of the type~$i$ such that particles of the $i$th type move along these contours, $i=1,\dots,N.$ The particles of the type $i$ move in the contours of the type~$i,$ $i=1,\dots,n.$ If $n=2,$ then there are $N_2$ contours 
of the type ~1, and there are $N_1$ contours of the type ~2. The $i$th contour of the type~1 contains  the cells $(0,i),(1,i),\dots,(N_1,i),$ $i=1,\dots,N_2.$
The $i$th contour of the type~2 contains  the cells $(i,0),(i,1),\dots,(i,N_2),$ $i=1,\dots,N_1.$    

In a Buslaev net, if more than one particle come to a common node simultaneously, then then a conflict occurs and a particle passes through the node according to the prescribed conflict resolution rule. In the considered net, a conflict occurs between a particles different type. In the contour networks, the conflict may be deterministic or stochastic. In the considered system a conflict is resolved by that first the particles of the first type move, and then the particles of the second type move, etc.

\section{
Sufficient condition for self-organization
}
\label{section:BuslN}
We say that, at time $t_0,$ the system is in a state of free movement if, at any moment $t=t_0,t_0+1,t_0+2,\dots,$ all particles move with probability~1.   

Let $d$ be the greatest common divisor of the numbers $N_1,$ 
$N_2,\dots,N_n.$

A set of cells $(x_1,\dots,x_n)$ such that 
$$\sum\limits_{i=1}^nx_i=k\ {\rm (mod}\ d)$$ 
is called the diagonal $k,$ $k=0,1,\dots,d-1.$

By definition, put
$$\phi_s(t)=k-t\ (mod\ d)$$
if, at time $t,$ the $s$th particle is in a cell of the diagonal  $k,$ $s=1,\dots,m,$ 
$k=1,\dots,d.$

Suppose $(x_1(t),\dots,x_d(t))$ is the vector such that $x_a(t)=0$ if, for any $s=1,\dots,m,$ $\phi_s(t)\ne a,$ and $x_k(t)=1$ if, there exists $s\in \{1,\dots,m\}$ such that $\phi_s(t)= a.$ 
\vskip 5pt
{\bf Lemma 1.} {\it Suppose $\phi_s(t)=a.$ If, at time $t,$ the $s$th particle 
moves, then $\phi_s(t+1)=a,$ and, if the particle does not move, then $\phi_s(t+1)=a-1$ (mod $d)$} 
\vskip 5pt
Lemma 1 follows from definitions.
\vskip 5pt
{\bf Lemma 2.} {\it Suppose $x_{a-1}(t)=1,$ $x_a(t)=x_{a+1}(t)=\dots=
x_{a+l}(t)=0$ (addition modulo $d),$ $l\ge 1.$ Then $x_{a-1}(t+1)=1,$ $x_a(t+1)=x_{a+1}(t+1),\dots,x_{a+l-1}(t+1)=0.$}
\vskip 5pt
{\bf Proof.} At time $t,$ the $s$th particle such that $\phi_s(t)=a-1$ moves. Therefore, $\phi_s(t+1)=a-1.$ Hence, 
$$x_{a-1}(t+1)=1.\eqno(1)$$ 
Suppose  there exists  $b\in \{0,1,\dots,l-1\}$ such that $x_{a+b}(t+1)=1.$ 
Then there exists $s\in \{1,\dots,m\}$ such that $\phi_s(t)=a+b$ or $\phi_s(t)=a+b+1.$
Hence there exists $c\in \{0,1,\dots,l\}$ such that $x_{a+c}(t)=1.$ This contradiction proves that 
$$x_a(t+1)=x_{a+1}(t+1),\dots,x_{a+l-1}(t+1)=0.\eqno(2)$$       
Combining (1), (2), we get  Lemma 2. 
\vskip 5pt
{\bf Lemma 3.} {\it The vector $X(t+1)$ contains not more 
clusters of zeros of length not less 2 in the vector than the vector $X(t).$ To any cluster of zeros of length not less than~2 in the vector $X (t+1)$  corresponds to a cluster of length not less than 2 with the same left limit in the vector $X (t)$.
}
\vskip 5pt
Lemma 3 follows from Lemmas 1 and 2.
\vskip 5pt
{\bf Lemma 4.} {\it If the system never results in the state of free movement, then, for any realization of the stochastic process, there exists a moment $t_0$ such that, for any $t\ge t_0,$ the vector $X(t)$ does not contain clusters of length less than~2.}
\vskip 5pt
{\bf Proof.} If the system never results in a free movement, then exists a particle (let it be a particle~$s)$ such that infinite number of delays occurs.  Therefore, as $t$ increases, the function $\phi_s(t)$ takes successively all the values $d-1,d-2,\dots,1,0.$ Suppose, after the end of this cycle, i.~e., at time such that a state is repeated, (suppose this cycle ends at time $T),$ in the vector $X(T)$ there are clusters of zeros of length not less 2. Therefore, taking into account Lemma~2, we get that, for any $t=0,1,\dots,T,$ the vector $X(t)$ contains a cluster of length not less than~2 with the same left limit. Hence the function $\phi_s(t)$ can not take all the values $d-1,d-2,\dots,1,0.$ This concludes the proof of the lemma.

\vskip 5pt
{\bf Theorem 1.} {\it If $m\le \frac{d}{2},$ then the system results in a state of free movement after a time interval with a finite expectation.}
\vskip 5pt
{\bf Proof.} If $m<d/2,$ then, for any $t,$ in the vector $X(t),$ there exists, at least one cluster of zeros such that there exists at least one cluster of length not less than~2. From this and Lemma~4, taking into account that the number of states is finite, we get the statement of the theorem under the assumption that $m<d/2.$

Suppose $m=d/2$ and the system does not result in the state of free movement. This condition may hold only if $d$ is an even number. According to Lemma 3 there exists a moment $t_0$ such that, for  $t\ge t_0,$ there exist only clusters of zeros of length~1 and clusters of zeros of length~1. If a particle is on the diagonal $k,$ then there are not other particles on the diagonal $k$ and particles on the diagonal $k+1.$  Hence, at any time $t\ge t_0,$ delays do not occur. This completes the proof of Theorem~1.

\section{Conclusion}

Versions of BML-models are helpful to obtain optimal solutions for problems of optimal traffic control.  

A system is studied such that, in this system, paticles move along an 
$n$-dimensional lattice  $N_1\times N_2\times \dots \times N_n$ with periodic boundary conditions according to rules generalizing the rules of the BML model. 
In Buslaev nets conception the self-organization property of the system corresponds to the existance of velocity single point spectrum equal to 1.
A sufficient condition for the self-organization of the system has been found. According to this condition, if number of particles is less than half of the greatest common divisor of the numbers $N_1$ $N_2,\dots,N_n$, then the system results in a state of free movement under the assumption that the number of particles is less than half the greatest common divisor.

\section*{Acknowledgments}

The authors would like to thanks
 our colleague Sergey  A. Buruykov for  the calculations
of the model realization.

\bibliographystyle{unsrt}  


\end{document}